\documentclass[runningheads]{llncs}
\usepackage[T1]{fontenc}

\usepackage{newtxtext}       %
\usepackage[varvw]{newtxmath}       %

\usepackage{mathtools}
\usepackage{xurl}
\usepackage{hyperref}
\usepackage{orcidlink}
\usepackage{booktabs}
\usepackage{wrapfig}
\usepackage{capt-of}

\newcommand{\R}{\mathbb{R}}
\newcommand{\N}{\mathbb{N}}
\newcommand{\Sph}{\mathbb{S}}
\newcommand{\Disk}{\mathbb{D}}
\newcommand{\Prob}{\mathcal{P}}
\newcommand{\Mes}{\mathcal{M}}
\newcommand{\Bor}{\mathcal{B}}
\newcommand{\Rad}{\mathcal{R}}
\newcommand{\NR}{\mathcal{N}}
\newcommand{\Wht}{\mathcal{W}}
\newcommand{\Met}{\mathfrak{M}}
\newcommand{\Gr}{\mathfrak{G}}
\newcommand{\Wass}{\operatorname{W}}
\newcommand{\SWass}{\operatorname{SW}}

\newcommand{\Four}{\mathcal{F}}
\newcommand{\Orbit}{\mathfrak{O}}
\newcommand{\dd}{\,\mathrm{d}}
\newcommand{\GL}{\mathrm{GL}}

\newcommand{\Leb}{\mathrm{L}}
\newcommand{\Orth}{\mathrm{O}}

\newcommand{\dist}{\operatorname{dist}}
\newcommand{\supp}{\operatorname{supp}}

\newcommand{\argmin}{\operatorname*{arg\,min}}

\newcommand{\mean}{\operatorname{mean}}
\newcommand{\std}{\operatorname{std}}

\newcommand{\T}{\operatorname{T}}
\newcommand{\Aff}{\operatorname{Aff}}
\newcommand{\op}{\operatorname{op}}
\newcommand{\aff}{\operatorname{aff}}

\newcommand{\Lebrho}{\scalebox{0.6}{\raisebox{-5pt}{$\Leb^2_\rho(\R)$}}}
\newcommand{\Lebrhou}{\scalebox{0.6}{\raisebox{-5pt}{$\Leb^2_{\rho\otimes u_{\Sph^1}}(\R\times\Sph^1)$}}}
\newcommand{\UniSph}[1]{\scalebox{0.6}{\raisebox{-2pt}{$\R\times\Sph^{#1}$}}}

\newcommand{\UniSphE}[1]{\R\times\Sph^{#1}}
\newcommand{\UniE}{\R}

\usepackage{graphicx}
\usepackage{color}

\newif\ifsidebysidelandscape
\sidebysidelandscapefalse
\usepackage{graphicx}

\begin{document}
\title{From Affine Transformations to Orthogonal Orbits: \\ A Radon–CDT Metric on Affine Classes} %

\titlerunning{From Affine Transformations to Orthogonal Orbits: R–CDT Metric on Affine Classes}

\author{Christos Sonis\inst{1,2}\orcidlink{0009-0007-6720-6555} \and
Jonas Bresch\inst{1,2}\orcidlink{0009-0009-3855-0933} \and
Gitta Kutyniok\inst{1,2,3}\orcidlink{0000-0001-9738-2487}}
\authorrunning{C. Sonis, J. Bresch, G. Kutyniok}
\institute{Ludwig-Maximilians-Universität München,
\href{https://www.math.lmu.de/math4ai/en/}{Department of Mathematics}
\email{\{sonis,bresch,kutyniok\}@math.lmu.de}
\and
\href{https://mcml.ai/research/groups/kutyniok/}{Munich Center for Machine Learning (MCML)},
\and
University of Tromsø
and
DLR-German Aerospace Center
}
\maketitle              %
\begin{abstract}
Normalized Radon cumulative distribution transforms (NR-CDTs)
provide affine-invariant image features,
but aggregation over projection directions
can assign identical features to distinct affine classes.
We ensure injectivity by whitening each measure and representing
its Radon--CDT through its orthogonal orbit.
For planar probability measures with finite second moments and
positive definite covariance, we prove that these orbits characterize
affine equivalence.
Optimal orthogonal alignment then defines a sliced Wasserstein
metric on affine classes, yielding separation margins and
perturbation guarantees for finite template families.
Numerical experiments demonstrate improved nearest-neighbor classification,
comparable accuracy with about half as many projection angles,
and accurate recovery of affine transformations.
\end{abstract}

\section{Introduction}

Optimal transport provides a geometric framework for comparing
probability measures through the displacement of mass
\cite{peyre2019computational,villani2003topics}.
In one dimension, the quadratic Wasserstein distance equals
the \(\Leb^2\)-distance between quantile functions.
This identity underlies the \emph{cumulative distribution transform} (CDT)
and, after integration over projection directions,
the sliced Wasserstein distance
\cite{bonneel2015sliced,park2018cdt}.
The \emph{Radon-CDT} (R-CDT)
uses this connection to represent images through
one-dimensional transport maps \cite{kolouri2016radon}.

For image comparison under affine distortions,
\emph{normalized R-CDT} (NR-CDT) variants provide invariant features
and classification guarantees with limited training data
\cite{beckmann2025max,beckmann2026generalizations,beckmann2026robustness}.
Their ability to distinguish affine classes, however, is limited.
We show that distinct affine classes can share the same set
of normalized slice values at every quantile level.
Consequently, every pointwise aggregation depending only
on these sets assigns identical features to the two classes.

We overcome this ambiguity by combining self-whitening with the R-CDT.
Centering and covariance normalization reduce invertible affine
transformations to orthogonal transformations
\cite{lan2013whitening}, which act on the projection directions
of the \emph{whitened R-CDT} (WR-CDT).
We prove that its orthogonal orbits characterize affine equivalence
for planar probability measures with finite second moments and
positive definite covariance.
Our measure-theoretic formulation includes atomic measures,
extending the original density-based R-CDT framework \cite{kolouri2016radon}.
Optimal orthogonal alignment defines a metric on affine classes,
equivalent to sliced Wasserstein alignment of the whitened measures,
and yields separation margins and perturbation guarantees
for finite template families.
Experiments demonstrate improved nearest-neighbor classification
over NR-CDT baselines, comparable accuracy with about half
as many projection angles in multiclass problems with limited training samples per class,
and accurate affine map recovery.

\section{Background and Limitation of the Normalized R-CDT}

We consider finite, signed, regular Borel measures $\Mes(\R^2)$ on $\R^2$.
The subset of probability measures are $\Prob(\R^2) \coloneqq \{\mu \in \Mes_+(\R^2) : \mu(\R^2) = 1\}$,
where $\Mes_+(\R^2) \subset \Mes(\R^2)$ is the subset of non-negative measures,
i.e. where $\mu(B) \geq 0$ for any Borel set $B \in \Bor(\R^2)$.
For a measurable map $T : \R^2 \to \R^2$, i.e. where $T^{-1}(B) \in \Bor(\R^2)$ for any $B \in \Bor(\R^2)$,
we denote the pushforward of $\mu \in \Mes(\R^2)$ on $T$ by $T_{\#}\mu \coloneqq \mu \circ T^{-1}$.
Hence, for any $T_{\#}\mu$-integrable function $g \in \Leb_{T_{\#}\mu}$, we have
$\int_{\R^2} g(x) \dd T_{\#} \mu(x) = \int_{\R^2} (g \circ T)(y) \dd \mu(y) < \infty$.
The probability measures with finite \emph{$p$th moment}
are defined by
\begin{equation*}
    \smash{\Prob\!\!_p(\R^2) \coloneqq \Bigl\{ \mu \in \Prob(\R^2)
                \;\Big|\;
                \scalebox{1.0}{$\int_{\R^2}$} \|x\|^p \dd \mu(x) < \infty\Bigr\},
    \quad p \geq 1,}
\end{equation*}
which include the probability measures $\Prob_{\!\mathrm c}(\R^2)$ with compact \emph{support},
where
\begin{equation*}
    \smash{\supp(\mu) \coloneqq \{x \in \R^2 \mid \mu(B_r(x)) > 0 \;\forall r > 0\},
    \quad B_r(x) \coloneqq \{y \in \R^2 \mid \|x - y\|_2 < r\}.}
\end{equation*}
By Hölder’s inequality, it holds $\Prob\!\!_p(\R^2)\subseteq\Prob\!\!_q(\R^2)$,
where $p\ge q\ge1$.
We further define the \emph{mean vector} and \emph{covariance matrix} of $\mu \in \Prob\!_p(\R^2)$ by
\begin{equation*}
    \smash{m_\mu \coloneqq \int_{\R^2}x\dd\mu(x),
    \qquad
    \Sigma_\mu \coloneqq \int_{\R^2}(x-m_\mu)(x-m_\mu)^\top\dd\mu(x),
    \quad p \geq 2.}
\end{equation*}

\subsection{Radon Projections and the Normalized R-CDT}

For $\theta\in\Sph^1$, let
$S_\theta:\R^2\to\R$ be given by
$S_\theta(x)\coloneqq\langle x,\theta\rangle$.
The \emph{restricted Radon transform} of $\mu\in\Prob(\R^2)$
in direction $\theta$ is the pushforward probability measure
$\Rad_\theta[\mu]\coloneqq(S_\theta)_\#\mu \in \Prob(\R)$.
For $\eta\in\Prob(\R)$,
its \emph{cumulative distribution function}
and \emph{quantile function} are
\begin{equation*}
    \smash{F_\eta(t)\coloneqq\eta((-\infty,t]),
    \quad t \in \R,
    \qquad
    F_\eta^{[-1]}(u)\coloneqq
    \inf\{t\in\R:F_\eta(t) > u\},
    \quad u\in(0,1).}
\end{equation*}
Fix an atomless reference measure $\rho\in\Prob(\R)$,
i.e. $\rho(\{s\})=0$ for every $s\in\R$.
The monotone rearrangement from $\rho$ to $\eta \in \Prob_2(\R)$
also-called the \emph{cumulative distribution transform} (CDT)
is given $\rho$-a.e.
by $Q_\eta \coloneqq F_\eta^{[-1]} \circ F_\rho \in \Leb^2_\rho(\R)$.
For any convex cost $c:\R\to[0,\infty)$,
provided $\int_{\R}c(s-Q_{\eta}(s))\,\mathrm d\rho(s) <\infty$,
the mapping $Q_{\eta}$ solves the Monge--Kantorovich \emph{optimal transport problem} \cite[§~2]{santambrogio2015optimal}:
\begin{equation}
    \label{eq: monge}
    \smash{Q_{\eta} \in
    \argmin_{
        \substack{T:\R\to\R\ \text{measurable}, %
        T_\#\rho=\eta}}
        \int_{\R}c(s-T(s))\dd \rho(s).}
\end{equation}
The monotone couplings \eqref{eq: monge} for $\eta,\xi \in \Prob_2(\R)$
provide the well-known Wasserstein-2 metric \cite{santambrogio2015optimal}
for the convex cost $c = |\cdot|^2$,
where $\Pi(\eta,\xi) \coloneqq \{\pi \in \Prob(\R^2) : \pi(\R\times B) = \xi(B), \pi(A\times \R) = \eta(A), A, B \in \Bor(\R)\}$:
\vspace{-10pt}
\begin{equation*}
    \Wass_2^2(\eta,\xi) \!
    \coloneqq \hspace{-9pt}\min_{\pi \in \Pi(\eta,\xi)} \int_{\R^2} \!\!\!c(x-y) \dd \pi(x,y)
    = \left\|Q_\eta - Q_\xi\right\|^2_{\Lebrho}
    \hspace{-8pt}\eqqcolon \!\!\int_\R\!\!|Q_\eta(u) - Q_\xi(u)|^2 \dd \rho(u).
    \vspace{-1pt}
\end{equation*}
For $\Rad_\theta[\mu] \in \Prob_2(\R)$,
the so-called \emph{Radon cumulative distribution transform} (R-CDT) \cite{kolouri2016radon}
is given $\rho$-a.e. by
\begin{equation}
    \label{eq: rcdt}
    \smash{Q_{\theta,\mu}
    \coloneqq
    Q_{\Rad_\theta[\mu]}
    = F_{\Rad_\theta[\mu]}^{[-1]}\circ F_\rho,
    \quad \text{with} \quad
    (Q_{\theta,\mu})_\#\rho=\Rad_\theta[\mu].}
\end{equation}
This immediately yields the correspondence for the \emph{sliced Wasserstein-2 distance} \cite{bonneel2015sliced}
\begin{align}
    \SWass_2^2(\mu, \nu)
    &\coloneqq \int_{\Sph^1} \Wass_2(\Rad_{\theta}[\mu], \Rad_{\theta}[\nu])^2 \dd u_{\Sph^1}(\theta) \label{eq: sliced wasserstein} \\
    &= \int_{\UniSph{1}} \hspace{-10pt}
    |Q_{\theta,\mu}(u) - Q_{\theta,\nu}(u)|^2 \dd(\rho\otimes u_{\Sph^1})(u,\theta)
    \eqqcolon \left\|Q_{\theta,\mu} - Q_{\theta,\nu}\right\|^2_{\Lebrhou}\!\!\!. \notag
\end{align}
The \emph{normalized R-CDT} (NR-CDT) \cite{beckmann2025max,beckmann2026generalizations} slices are defined as
\begin{equation}
    \label{eq: nrcdt}
    \smash{\NR_\theta[\mu](u)
    \coloneqq
    \tfrac{Q_{\theta,\mu}(u)-m_{\theta,\mu}}
    {\sigma_{\theta,\mu}},
    \qquad u \in \R,
    \qquad \theta \in \Sph^1,}
\end{equation}
where $m_{\theta,\mu} \coloneqq \int_{\R} Q_{\theta,\mu}(u)\dd \rho(u)$
and $\sigma_{\theta,\mu}^2 \coloneqq \int_{\R} |Q_{\theta,\mu}(u)-m_{\theta,\mu}|^2\dd \rho(u)$.
The NR-CDT projections \eqref{eq: nrcdt}
are well-defined \cite[Prop.~4]{beckmann2026generalizations} for any measure in
\begin{equation}
    \label{eq: porb c *}
    \smash{\Prob_{\!\mathrm c}^*(\R^2)
    \coloneqq \bigl\{\mu \in \Prob_{\mathrm c}(\R^2)
    \;\big|\;
    \sigma_{\theta,\mu}^2 > 0
    \;\; \forall \theta \in \Sph^1\bigr\}.}
\end{equation}

\begin{lemma}
    \label{lem: mean,std}
    For $\mu \in \Prob_{2}(\R^2)$
    and $\theta \in \Sph^1$,
    it holds $m_{\theta,\mu} = \mean(Q_{\theta, \mu}) = \langle m_\mu,\theta\rangle$
    and $\sigma_{\theta,\mu}^2 = \std(Q_{\theta, \mu})^2 = \theta^\top\Sigma_\mu\theta$.
\end{lemma}
\begin{proof}
    The first equality in each identity is displayed in the proof of \cite[Lem.~4.2]{beckmann2026robustness};
    the second follows from the linearity of the integral.
    The second statement follows by $(Q_{\theta,\mu})_{\#}\rho = \Rad_{\theta}[\mu]$
    and
    \begin{align*}
        \smash{\theta^\top \Sigma_\mu \theta
        = \!\!\int_{\R^2} \bigl(\langle x, \theta\rangle - \langle m_\mu, \theta \rangle\bigr)^2 \dd \mu(x)
        = \!\!\int_{\R} (s - m_{\theta,\mu})^2 \dd \overbracket{(S_\theta)_{\#} \mu}^{= \Rad_\theta[\mu]}(s)
        = \std(Q_{\theta,\mu})^{2}.}
        \hspace{-10pt}\tag*{\qed}
    \end{align*}
\end{proof}

\begin{proposition}
    \label{prop: sigma support}
    The following three statements are equivalent for $\mu \in \Prob_2(\R^2)$: \\[-12pt]
    \begin{center}
        \textup{i)} $\exists c > 0$ s.t. $\sigma_{\theta,\mu}^2 \geq c$ for all $\theta \in \Sph^1$,
        \quad
        \textup{ii)} $\dim(\aff(\supp(\mu))) > 1$,
        \quad
        \textup{iii)} $\Sigma_\mu \succ 0$.
    \end{center}
\end{proposition}
\begin{proof}
    By Lem.~\ref{lem: mean,std}, $\sigma_{\theta,\mu}^2=\theta^\top\Sigma_\mu\theta\geq\lambda_{\min}(\Sigma_\mu)$ with equality at a unit eigenvector, so (i) and (iii) are equivalent.
    For the equivalence of (ii) and (iii),
    we first note that $\Sigma_\mu \succeq0$ by construction.
    If we assume that $\Sigma_\mu \not\succ0$,
    then there exists some $y \in \R^2\setminus\{0\}$ such that $y^\top \Sigma_\mu y = 0$.
    Hence, we would have
    \begin{equation}
        \label{eq: one dimensional sigma}
        \smash{0 = y^\top \Sigma_\mu y
        = \int_{\R^2} y^\top (x - m_\mu)(x - m_\mu)^\top y \dd \mu(x)
        = \int_{\R^2} (y^\top (x - m_\mu))^2 \dd \mu(x).}
    \end{equation}
    Therefore, we need to have $\supp(\mu) \subseteq \{x \in \R^2 : y^\top(x - m_\mu) = 0\}$,
    since this line is closed and has full $\mu$-measure.
    Conversely, if we assume that $\dim(\aff(\supp(\mu))) \leq 1$,
    then $\supp(\mu) \subseteq \{x \in \R^2 : y^\top x = b\}$
    for some $y \in \R^2 \setminus\{0\}$, $b\in\R$.
    Since $\mu$ is concentrated on this line,
    $y^\top m_\mu = b$.
    By \eqref{eq: one dimensional sigma},
    it follows that $y^\top \Sigma_\mu y = 0$.
    \qed
\end{proof}

By Prop.~\ref{prop: sigma support},
we obtain for the considered set of probability measures \eqref{eq: porb c *}
that
\begin{align}
    \label{eq: measure 2*}
    \smash{\Prob_{\rm c}^*(\R^2)
    = \{\mu \in \Prob_{\rm c}(\R^2) \mid \Sigma_{\mu} \succ 0\}
    \subset \Prob_2^*(\R^2)
    \coloneqq \{\mu \in \Prob_{2}(\R^2) \mid \Sigma_\mu \succ 0\}.}
\end{align}
Hence,
we enlarge the domain by whitening \eqref{eq: whitening},
including those measures with noncompact support.

\subsection{Affine Measure Classes}

In contrast to \cite{beckmann2025max,beckmann2026generalizations,beckmann2026robustness},
we specifically treat the measure classes
as \emph{orbits} of the \emph{affine group}
$\smash{\Aff_2(\R)
=
\left\{
\alpha_{A,y} : \R^2 \to \R^2
\mid \alpha_{A,y}=A\cdot+y,\;
A\in \GL_2(\R),\;
y\in\R^2
\right\}}$,
where $\GL_2(\R) \coloneqq \{A \in \R^{2\times2} \mid \exists A^{-1} \in \R^{2\times2}, AA^{-1} = A^{-1}A = I_2\}$,
operating by push forward
$\alpha_{A,y}(\mu)
\coloneqq (\alpha_{A,y})_{\#}\mu
= \mu \circ \alpha_{A,y}^{-1}
= \mu \circ \bigl(A^{-1}(\cdot - y)\bigr)$,
where $\alpha_{A,y} \in \Aff_2(\R)$.

\begin{lemma}
    \label{lem: aff mean sigma}
    Let \(\mu \in \Prob_2(\mathbb R^2)\) and $\nu = (\alpha_{A,y})_{\#}\mu$, for some $\alpha_{A,y} \in \Aff_2(\R)$.
    Then $m_\nu = A m_\mu + y$ and $\Sigma_\nu = A \Sigma_\mu A^\top$.
\end{lemma}
\begin{proof}
    Direct calculations give
    $m_\nu = m_{(A\,\cdot\, + y)_{\#}\mu} = \int_{\R^2} A x + y \dd \mu(x) = A m_\mu + y$
    and
    \begin{align*}
        \hspace{-10pt}
        \Sigma_\nu
        = \Sigma_{(A\,\cdot\, + y)_{\#}\mu}
        = \!\!\int_{\R^2}\!\! (A x + y - m_{(A\,\cdot\, + y)_{\#}\mu})(A x + y - m_{(A\,\cdot\, + y)_{\#}\mu})^\top \!\dd \mu(x)
        = A \Sigma_\mu A^\top\!\!.
        \hspace{-10pt}
        \tag*{\qed}
    \end{align*}
\end{proof}

\begin{proposition}
    \label{prop: well}
    The affine group $\Aff_2(\R)$
    operates on $\Prob^*_2(\R^2)$
    by push forward %
    and $\Prob_2^*(\R^2)$ is closed under this operation.
\end{proposition}
\begin{proof}
    By \eqref{eq: measure 2*}, %
    \smash{$\Prob^*_2(\R^2) =\left\{\mu\in\mathcal{P}_2(\R^2):\Sigma_\mu\succ 0\right\}$}.
    Let $\alpha, \beta \in \Aff_2(\R)$, then
    \begin{equation*}
        \smash{(\alpha\circ\beta)_\#\mu =\mu\bigl((\alpha\circ\beta)^{-1}\bigr) =\mu\bigl(\beta^{-1}(\alpha^{-1}\bigr) =\beta_\#\mu\bigl(\alpha^{-1}\bigr)=\alpha_\#(\beta_\#\mu),
        \quad \operatorname{id}_{\#}\mu=\mu}
    \end{equation*}
    The operation is well defined,
    since affine images preserve finite 2nd moments,
    and with Lem.~\ref{lem: aff mean sigma}
    holds $\Sigma_{\alpha_{A,y}(\mu)} = A\Sigma_\mu A^\top \succ 0$
    for all $A \in \GL_2(\R)$
    if and only if $\Sigma_\mu \succ 0$.
    The latter follows since $x \mapsto A x$ is a bijection,
    and hence $x^\top A S A^\top x = (A^\top x)^\top S A^\top x > 0$
    for all $x\in\R^2\setminus\{0\}$ if and only if $S \succ 0$.
    \qed
\end{proof}

The affine class of $\mu$ is the orbit of the affine group %
by Prop.~\ref{prop: well}
\begin{equation}
    \label{eq: aff class}
    \smash{[\mu]
    \coloneqq
    \Aff_2(\R)(\mu)
    = \left\{(\alpha_{A,y})_{\#}\mu \coloneqq (A\cdot+y)_\#\mu \;\middle|\; A \in \GL_2(\R),\ y\in\R^2 \right\},}
\end{equation}
so $\nu \in [\mu]$, i.e. $\nu \sim \mu$,
is an equivalence relation on $\Prob^*_2(\R^2)$.
The set of affine classes can be expressed as
\begin{equation}
    \label{eq: aff class fac}
    \smash{\mathcal{A}\coloneqq \Prob_2^*(\R^2)/\Aff_2(\R).}
\end{equation}
Using the quantities $a_A(\theta) \coloneqq \|A^\top\theta\|$,
$b_y(\theta) \coloneqq \langle y,\theta\rangle$,
and $\phi_A(\theta) \coloneqq \tfrac{A^\top\theta}{\|A^\top\theta\|}$,
for some $A \in \GL_2(\R)$ and $y \in \R^2$,
the affine covariance identities, read as \cite[Prop.~3 \& 5]{beckmann2026generalizations}
\\[-15pt]
\begin{align}
    \Rad_\theta[(\alpha_{A,y})_{\#}\mu]
    &=
    (a_A(\theta)\,\cdot+b_y(\theta))_\#
    \Rad_{\phi_A(\theta)}[\mu]
    =
    \Rad_{\phi_A(\theta)}[\mu]\bigl(\tfrac{\;\cdot\; - b_y(\theta)}{a_A(\theta)}\bigr),
    \label{eq: radon-covariance}\\[-3pt]
    \NR_\theta[(\alpha_{A,y})_{\#}\mu]
    &=
    \NR_{\phi_A(\theta)}[\mu].
    \label{eq: nrcdt-covariance}
\end{align}
\vspace{-20pt}

\subsection{Non-injectivity of Normalized R-CDTs}
\label{sec: non-inj nrcdt}

We define for $\mu \in \Prob_{\rm c}^*(\R^2)$ the \emph{$h$-Normalized R-CDT} ($_h$NR-CDT) as in~\cite[§~3.4]{beckmann2026generalizations}:
\begin{equation}
    \label{eq: hnrcdt}
    \smash{\NR_h[\mu]\colon\R\to\R,
    \qquad \NR_h[\mu](u)\coloneqq h\big(\NR[\mu](u,\cdot)\big),
    \qquad u\in\R.}
\end{equation}
If $\eta:\R\to\R$ is a boundedness-preserving function and $H:\mathcal{B}(\R)\to\R$ is bounded,
where $\mathcal{B}(\R)$ denotes the set of bounded subsets of $\R$,
let \smash{$h(g)=H\left(\left\{\eta(g(\theta))\mid \theta\in \Sph^1\right\}\right)$}
for any $g\in L^\infty(\Sph^1)$
so that $h(g\circ\theta)=h(g)$
is automatically satisfied for all $g\in L^\infty(\Sph^1)$
and all bijections $\theta:\Sph^1\to \Sph^1$.
Therefore, the $_h$NR-CDTs \eqref{eq: hnrcdt} become
\begin{equation*}
    \smash{\NR_h[\mu](u) = H(\eta \left(I_\mu(u)\right)),
    \qquad \text{where} \qquad
    I_\mu(u) \coloneqq \left\{\NR_\theta[\mu](u) \mid \theta\in \Sph^1\right\}.}
\end{equation*}
For $\mu \in \Prob_{\rm c}^*(\R^2)$,
the $_{\mathrm{m}}$NR-CDT of~\cite[§~3.3]{beckmann2026generalizations} is obtained by $\eta = \operatorname{Id}$ and $H = \sup$.

\begin{figure}[t]
    \includegraphics[width=0.22\linewidth, trim={38pt 352pt 1005pt 0}, clip]{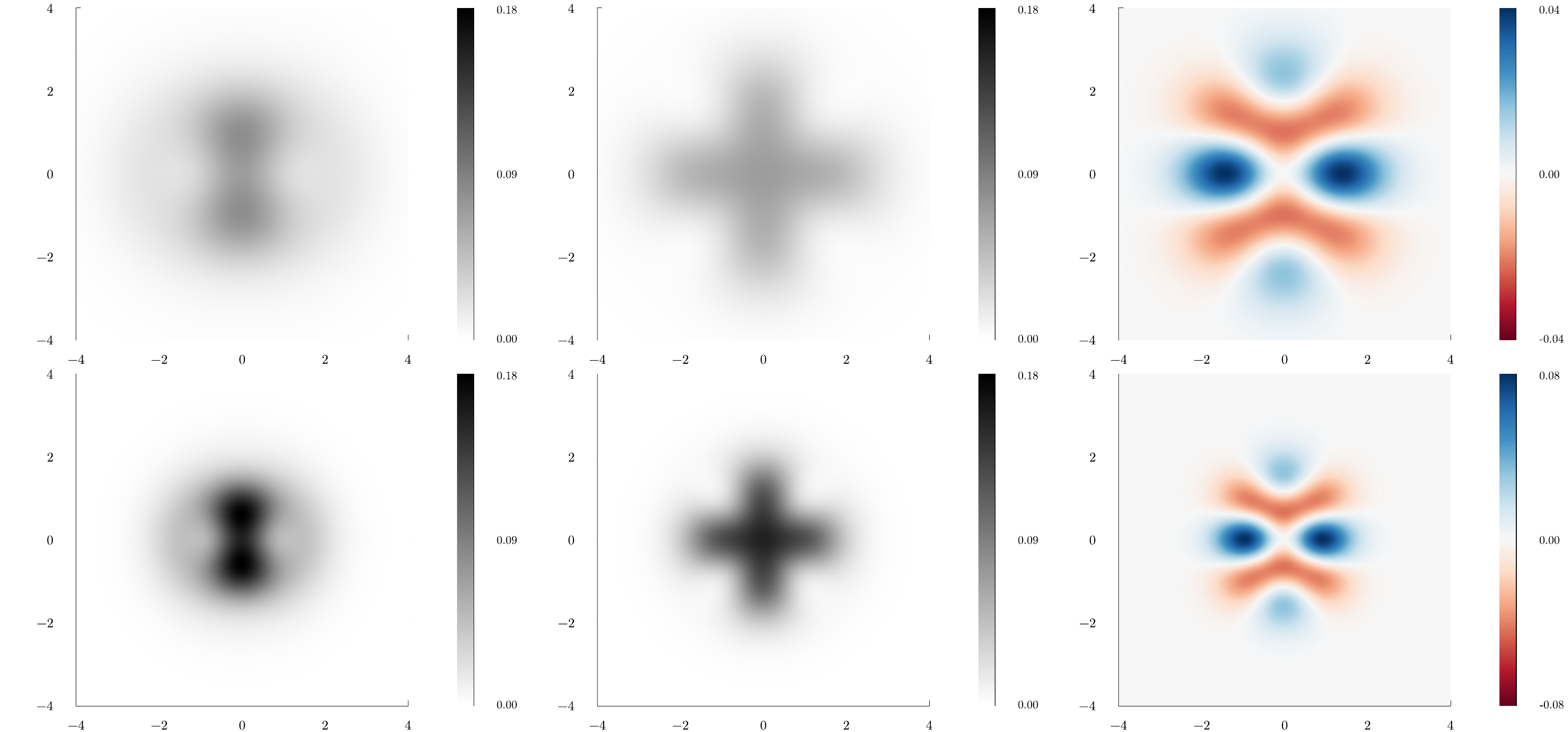}%
    \;
    \includegraphics[width=0.22\linewidth, trim={38pt 0pt 1005pt 350pt}, clip]{numeric/densities.pdf}%
    \;
    \includegraphics[width=0.25\linewidth, trim={4pt 670pt 1000pt 0}, clip]{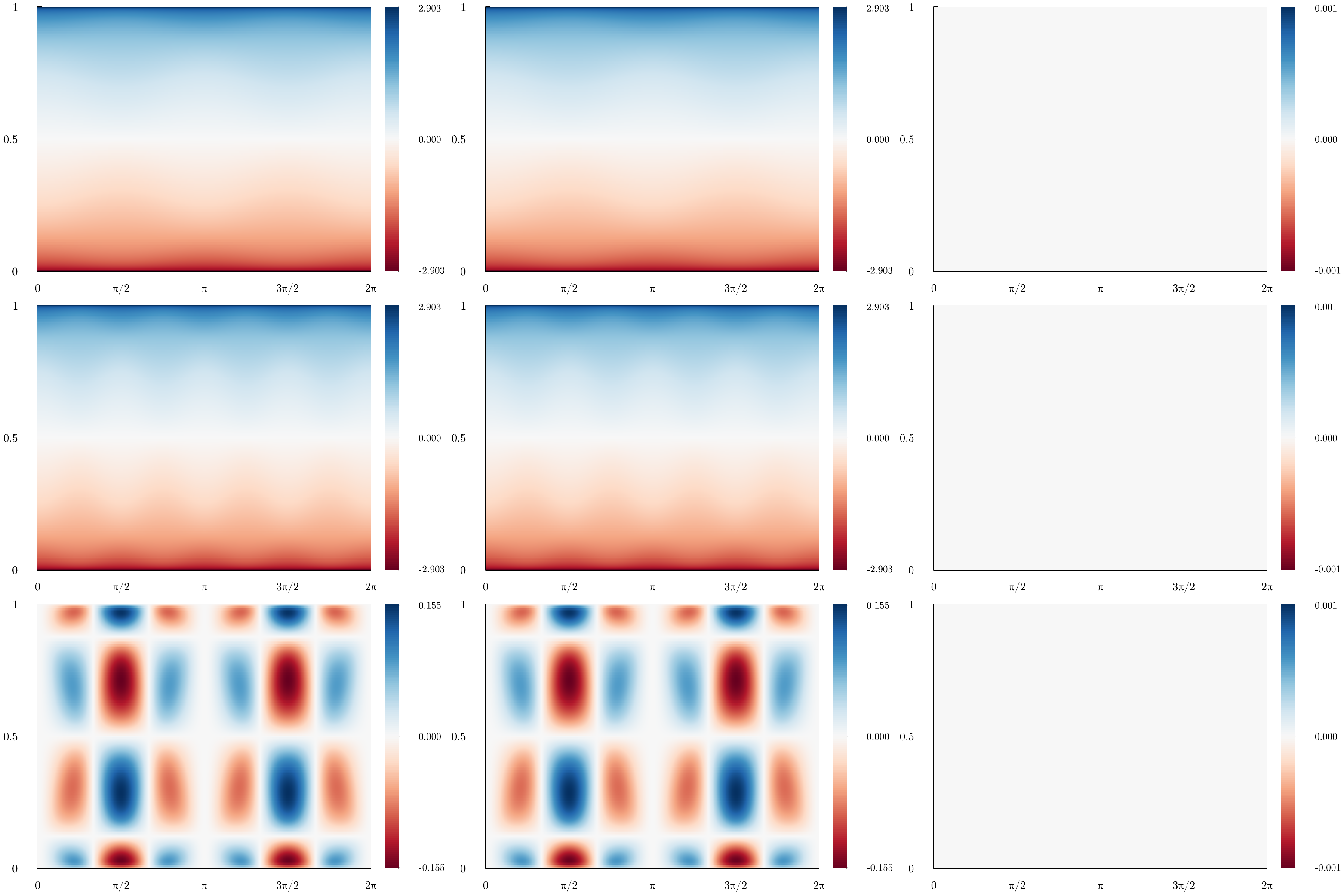}%
    \;
    \includegraphics[width=0.27\linewidth, trim={35pt 35pt 28pt 42pt}, clip]{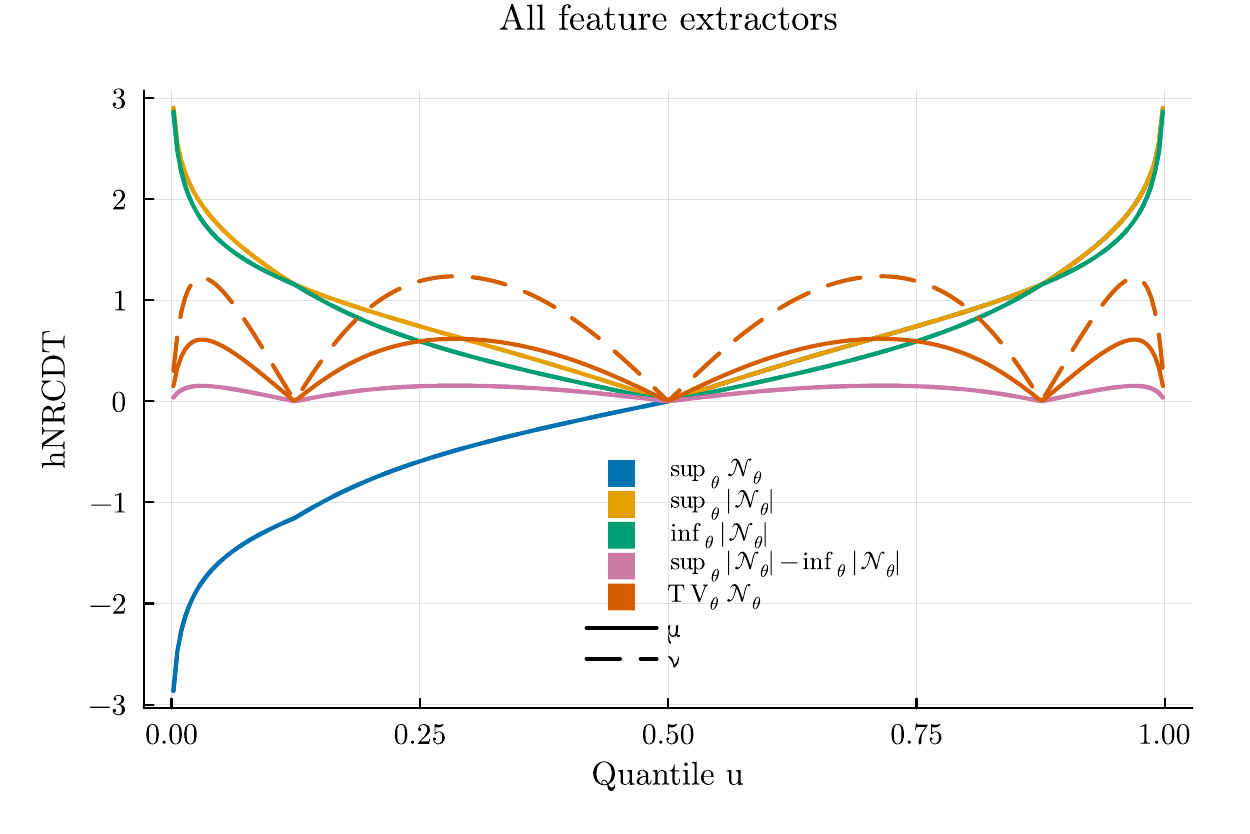}
    \includegraphics[width=0.22\linewidth, trim={538pt 350pt 505pt 0}, clip]{numeric/densities.pdf}%
    \;
    \includegraphics[width=0.22\linewidth, trim={538pt 0pt 505pt 350pt}, clip]{numeric/densities.pdf}%
    \;
    \includegraphics[width=0.25\linewidth, trim={5pt 330pt 1000pt 330pt}, clip]{numeric/nrcdts.pdf}%
    \;
    \includegraphics[width=0.27\linewidth, trim={65pt 40pt 25pt 42pt}, clip]{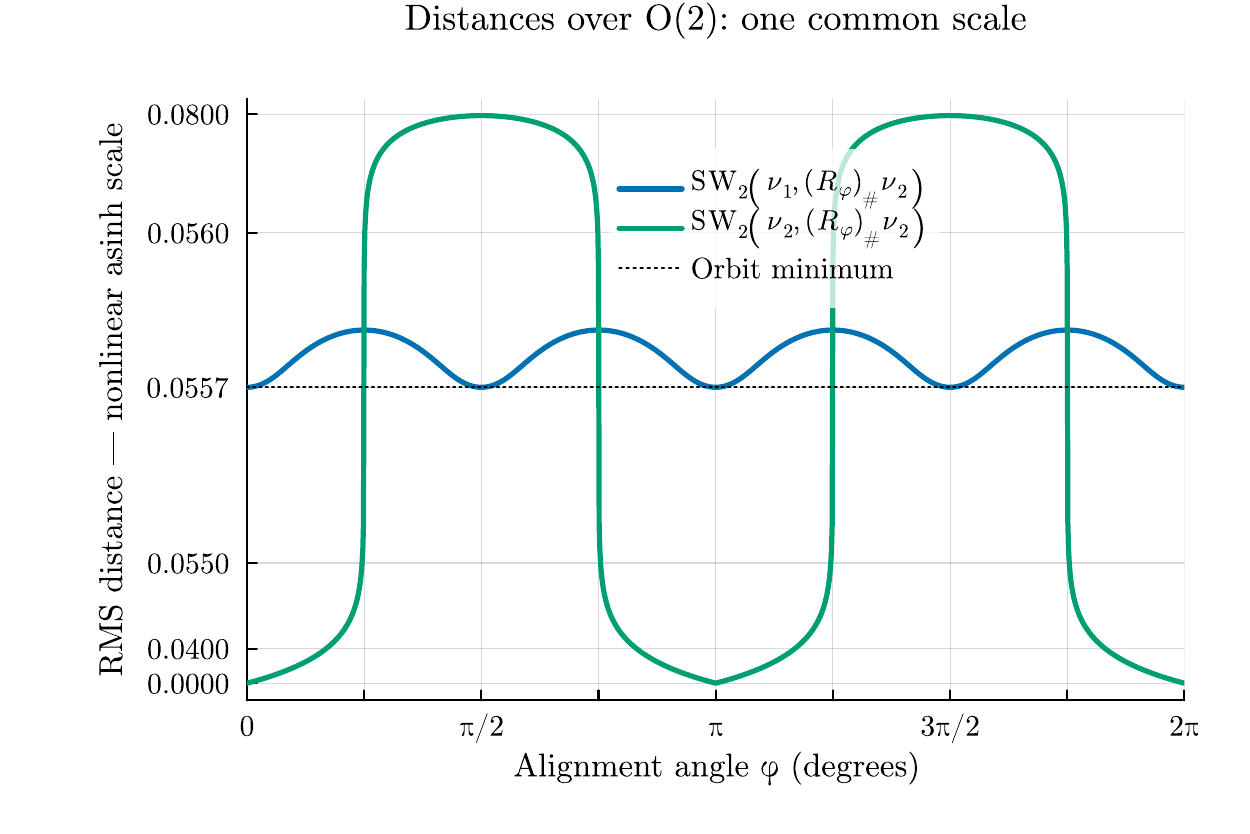} \\[-20pt]
    \caption{Template measures $\mu, \nu$ (\emph{left}),
    the whitened versions $\Wht[\mu], \Wht[\nu]$ (\emph{middle left}),
    the NR-CDTs $\NR[\mu],\NR[\nu]$ (\emph{middle right}).
    The $_h$NRCDTs
    and the distance in the orbits (\emph{right top/bottom}).}
    \label{fig: templates non-injective}
    \vspace{-10pt}
\end{figure}

\begin{example}
    Let $0 < a < b < c < d$ with $a + b + c + d = 2\pi/3$.
    Further,
    let $\mu$ be the uniform measure on the points $x_j \coloneqq (\cos\beta_j, \sin\beta_j)$, $j \in \{1, \dots, 12\}$,
    with $\beta_1 = 0$, $\beta_{13} = 2\pi$ and consecutive gaps $\beta_{j+1} - \beta_j$ equal to $(a,b,c,d)$ repeated three times.
    Additionally, let $\nu$ be defined in the same way
    with gaps $(a,b,d,c)$ repeated three times.
    One can show that $\mu, \nu \in \Prob_{\rm c}^*(\R^2)$,
    $\nu \notin [\mu]$, and $\NR_{\mathrm m}[\mu] = \NR_{\mathrm m}[\nu]$ \(\rho\)-a.e. on \(\mathbb R\).
\end{example}

Therefore, the $_\textrm{m}$NR-CDT, which is constant on affine classes \cite[Prop.~7]{beckmann2026generalizations},
induces a map $[\mu] \mapsto \NR_{\mathrm m}[\mu]$ on the set of affine classes $\Prob_{\rm c}^*(\R^2)/\Aff_2(\R)$
that is not injective.
Furthermore, by showing the continuity of the NR-CDT,
any $_h$NR-CDT \eqref{eq: hnrcdt} is not injective.
Fig.~\ref{fig: templates non-injective} illustrates the ambiguity for continuous inputs,
except for $_{\mathrm{tv}}$NR-CDT.

\section{Self-Whitening and the Affine-Canonical Transform}

\begin{figure}[t]
    \resizebox{\linewidth}{!}{%
    \includegraphics[width=0.2\linewidth,trim={0pt 305pt 815pt 0pt}, clip]{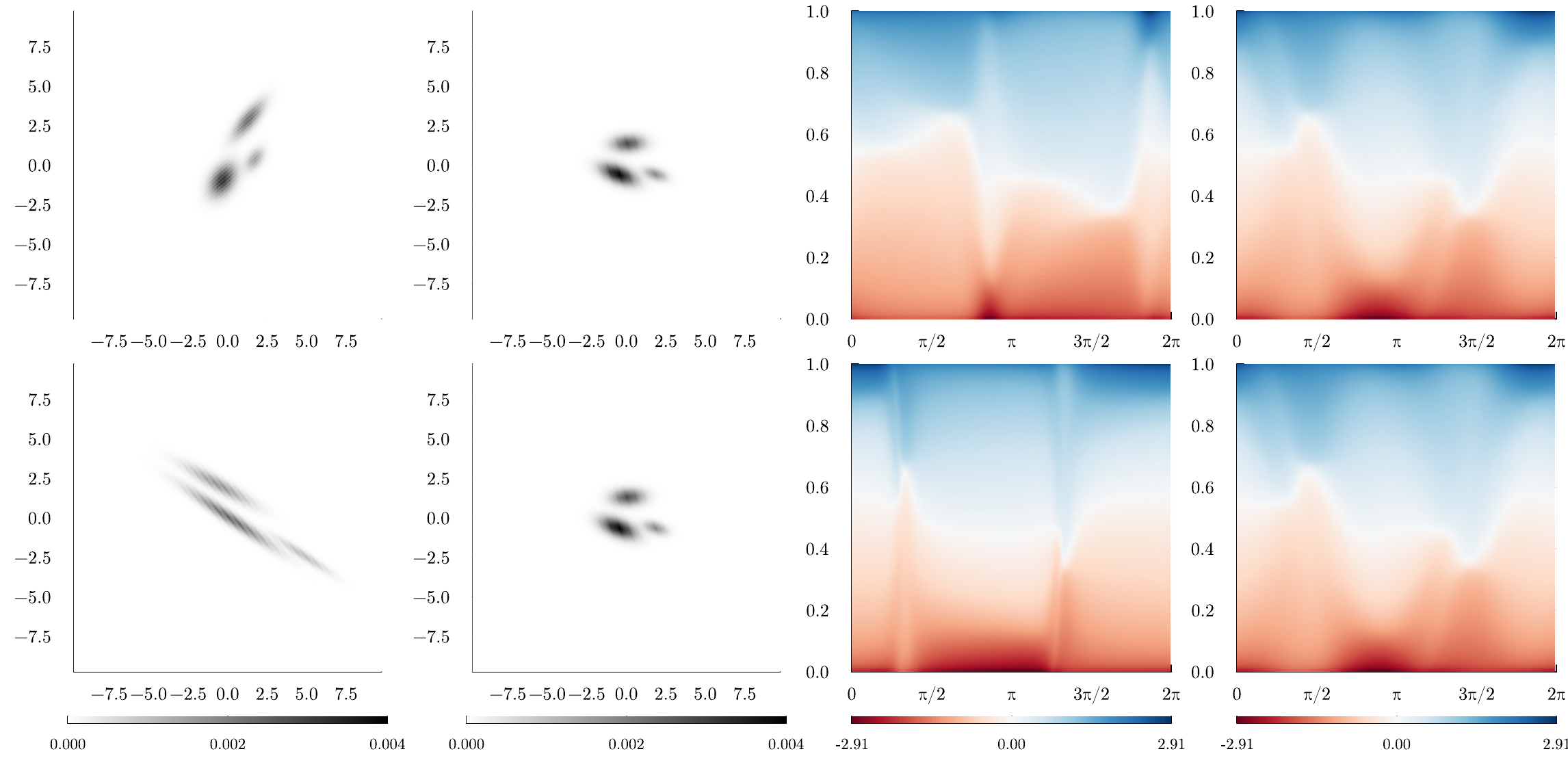}%
    \includegraphics[width=0.17\linewidth,trim={40pt 55pt 815pt 250pt}, clip]{numeric/whitening_desnities.pdf}%
    \includegraphics[width=0.16\linewidth,trim={320pt 305pt 555pt 0pt}, clip]{numeric/whitening_desnities.pdf}%
    \rotatebox{90}{\includegraphics[width=0.18\linewidth,trim={25pt 10pt 810pt 500pt}, clip]{numeric/whitening_desnities.pdf}}%
    \includegraphics[width=0.2\linewidth,trim={555pt 305pt 265pt 0pt}, clip]{numeric/whitening_desnities.pdf}%
    \includegraphics[width=0.178\linewidth,trim={585pt 56pt 265pt 250pt}, clip]{numeric/whitening_desnities.pdf}%
    \includegraphics[width=0.175\linewidth,trim={855pt 305pt 0pt 0pt}, clip]{numeric/whitening_desnities.pdf}%
    \rotatebox{90}{\;\includegraphics[width=0.17\linewidth,trim={585pt -10pt 265pt 500pt}, clip]{numeric/whitening_desnities.pdf}}%
    }\\[-20pt]
    \caption{Affinely transformed Gaußian mixtures (\emph{Left block: left and middle}),
    their whitened versions (\emph{Left block: right}),
    NR-CDT (\emph{Right block: left and middle}),
    and WR-CDT (\emph{Right block: right}).}
    \label{fig: whitening_densities}
    \vspace{-15pt}
\end{figure}

For $\mu \in \mathcal{P}_2^*(\R^2)$,
let $\Sigma_\mu^{1/2}$ be the symmetric positive definite square root of $\Sigma_\mu$,
well-defined by Prop.~\ref{prop: well},
and define the \emph{whitening map} and the \emph{whitened measure}
as
\begin{equation}
    \label{eq: whitening}
    \smash{\Wht[\mu] \coloneqq (\T_\mu)_\#\mu,
    \qquad \text{where} \qquad
    \T_\mu \coloneqq \Sigma_\mu^{-1/2}(\,\cdot-m_\mu).}
\end{equation}
Consequently,
a probability measure is called \emph{whitened},
i.e. standard normalized, if $m_\mu = 0$
and $\Sigma_{\mu} = I_2$.
In Fig.~\ref{fig: whitening_densities},
the action of the WR-CDT and the consequence for the NR-CDTs
are visualized.

In the following we mostly consider orthogonal group action %
i.e. where $(\alpha_{O,0})_{\#}\mu = O_{\#}\mu$
for some $O \in \Orth_2(\R) \coloneqq \{A \in \R^{2\times2} \mid A^\top A = A A^\top = I_2\} \subset \GL_2(\R)$.

\begin{lemma}
    \label{lem: orthogonal whitened}
    Let $\mu \in \mathcal{P}_2^*(\R^2)$.
    Then,
    $\Wht[\mu] \in \Prob_2^*(\R^2)$ is whitened, i.e. $m_{\Wht[\mu]} = 0$ and $\Sigma_{\Wht[\mu]} = I_2$.
    Moreover,
    $O_\# \mu$ is whitened whenever $\mu$ is whitened where $O \in \Orth_2(\R)$.
\end{lemma}
\begin{proof}
    The whitened measure $\Wht[\mu]$ is still in $\Mes_+(\R^2)$.
    The map $\T_\mu$ is the affine bijection \smash{$x \mapsto \Sigma_\mu^{-1/2}x-\Sigma_\mu^{-1/2}m_\mu$},
    and hence $\Wht[\mu] \in \Prob(\R^2)$ as well as by the integral transform w.r.t. $\T_\mu$
    it is of finite $2$nd moment.
    By Lem.~\ref{lem: aff mean sigma},
    \smash{$m_{\Wht[\mu]} = \Sigma_\mu^{-1/2}m_\mu - \Sigma_\mu^{-1/2}m_\mu = 0$}
    and $\smash{\Sigma_{\Wht[\mu]} = \Sigma_\mu^{-1/2}\Sigma_\mu\Sigma_\mu^{-1/2} = I_2}$,
    using the symmetry of \smash{$\Sigma_\mu$} and hence of \smash{$\Sigma_\mu^{1/2}$},
    which remains by Prop.~\ref{prop: well} positive definite%
    ---yielding the first assertion.
    For the second assertion,
    let $O \in \Orth_2(\R)$.
    By Lem.~\ref{lem: aff mean sigma},
    we obtain $m_{O_\#\mu} = O m_\mu = 0$,
    if and only if $m_\mu = O^{-1} 0 = 0$,
    which requires at least $O \in \GL_2(\R)$ and
    as stated $O \in \Orth_2(\R)$ for
    \begin{equation*}
        \smash{\Sigma_{O_\#\mu} \overset{\textrm{Lem.~\ref{lem: aff mean sigma}}}{=} O \Sigma_\mu O^\top = I_2
        \Leftrightarrow
        \Sigma_\mu = O^{-1} O^{-\top} = (O^\top O)^{-1} = I_2.}
        \tag*{\qed}
    \end{equation*}
\end{proof}

\begin{proposition}
    \label{prop: whitening}
    Let $\mu, \nu \in \Prob_2^*(\R^2)$
    and $\nu = (\alpha_{A,y})_{\#}\mu \in [\mu]$,
    for some $\alpha_{A,y} \in \Aff_2(\R)$.
    Then $\Wht[\nu] = O_\#\Wht[\mu]$,
    where $O = \Sigma_{\nu}^{-1/2}A\Sigma_{\mu}^{1/2}\in \Orth_2(\R)$.
\end{proposition}
\begin{proof}
    By Lem.~\ref{lem: aff mean sigma},
    it holds $m_{\nu} = A m_\mu + y$
    for some $A \in \GL_2(\R)$ and $y \in \R^2$ such that
    \begin{align*}
        \smash{T_{\nu} \circ (A\cdot + y)
        = \Sigma_{\nu}^{-1/2}(A\cdot + y - m_{\nu})
        = \Sigma_{\nu}^{-1/2}(A (\;\cdot\; - m_{\mu}))
        = O \circ \T_\mu.}
    \end{align*}
    To this end, with \eqref{eq: whitening},
    holds $\Wht[\nu] = (T_{\nu} \circ (A \cdot + y))_{\#} \mu = (O \circ \T_\mu)_{\#} \mu = O_{\#} \Wht[\mu].$
    By Lem.~\ref{lem: orthogonal whitened}, both $\Wht[\mu]$ and $\Wht[\nu]$ are whitened measures.
    By Lem.~\ref{lem: aff mean sigma} follows
    $I_2 = \Sigma_{\Wht[\nu]} = O \Sigma_{\Wht[\mu]} O^\top = O O^\top$.
    \qed
\end{proof}

\begin{proposition}
    \label{prop: whitened to affine}
    Let $\mu, \nu \in \mathcal{P}_2^*(\R^2)$,
    $\Wht[\nu] = O_\# \Wht[\mu]$ for $O \in \Orth_2(\R)$.
    Then $\nu \in [\mu]$.
\end{proposition}
\begin{proof}
    Since $T_\nu$ is a bijection,
    and since $S_\#(T_\#\mu) = (S \circ T)_\#\mu$,
    it follows that
    \begin{equation*}
        \smash{\nu
        = (T_\nu^{-1})_\#\Wht[\nu]
        = (T_\nu^{-1})_\# O_\#\Wht[\mu]
        = (T_\nu^{-1})_\# O_\# (\T_\mu)_\#\mu
        = (T_\nu^{-1}\circ O\circ T_\mu)_\#\mu,}
    \end{equation*}
    where %
    $T_\nu^{-1} \circ O \circ \T_\mu = \Sigma_\nu^{1/2}O\Sigma_\mu^{-1/2}(\,\cdot\, -m_\mu)+m_\nu = A\cdot+y$,
    using $A\coloneqq\Sigma_\nu^{1/2}O\Sigma_\mu^{-1/2}\in \GL_2(\R)$
    and $y \coloneqq m_\nu-Am_\mu$.
    Thus, $\nu = (A\cdot + y)_\#\mu
    \in[\mu]$.
    \qed
\end{proof}

From the latter statement
we can immediately derive for a canonical \emph{whitened direction} $\theta\in\Sph^1$
the corresponding direction in the original measure $\mu \in \Prob_2^*(\R^2)$.

\begin{corollary}
    \label{cor: no-resampling}
    For every $\mu, \nu\in \Prob_2^*(\R^2)$
    where $\nu = (\alpha_{A,y})_{\#}\mu \in [\mu]$ and $\theta\in\Sph^1$,
    it holds \\[-20pt]
    \begin{align*}
        \label{eq: no-resampling}
        \NR_\theta[\Wht[\mu]]
        &=
        \NR_{\theta_\mu(\theta)}[\mu],
        \qquad\qquad\qquad\qquad
        \theta_\mu(\theta)
        \coloneqq \tfrac{\Sigma_\mu^{-1/2}\theta}{\|\Sigma_\mu^{-1/2}\theta\|},
        \\[-10pt]
        \NR_\theta[\Wht[\nu]]
        &=
        \NR_{\theta_\mu(O^\top\theta)}[\mu]
        =
        \NR_{\eta_{(\alpha_{A,y})_{\#}\mu}(\theta)}[\mu],
        \qquad
        \eta_{(\alpha_{A,y})_{\#}\mu}(\theta)
        \coloneqq \tfrac{\Sigma_\mu^{-1} A^{-1} \Sigma_{\nu}^{1/2} \theta}{\|\Sigma_\mu^{-1} A^{-1} \Sigma_{\nu}^{1/2} \theta\|_2}.
    \end{align*}
\end{corollary}

According to Cor.~\ref{cor: no-resampling},
we define the \emph{Whitened-Radon-Cumulative-Distribution-Trans\-form} (WR-CDT)
for some measure $\mu \in \Prob_2^*(\R^2)$ from \eqref{eq: measure 2*} by
\begin{equation}
    \label{eq: wrcdt}
    \smash{X_\mu(\cdot, \theta)
    =\NR_{\theta}[\Wht[\mu]]
    \in
    \Leb^2_\rho(\UniE),
    \qquad
    \theta \in \Sph^1}
\end{equation}
and hence $X_{\mu} \in \Leb^2_{\rho \otimes u_{\Sph^1}}(\UniSphE{1})$.
Applying Cor.~\ref{cor: no-resampling},
we have the following correspondences
for the affine measure classes \eqref{eq: aff class}:
\begin{equation}
    \label{eq: o2-equivariance}
    \smash{X_\nu(\theta) = X_\mu(O^\top\theta),
    \qquad \nu \in [\mu]
    \quad \text{for some} \quad
    O \in \Orth_2(\R).}
\end{equation}
Moreover, by Lem.~\ref{lem: mean,std},
we have for whitened \smash{$\mu \in \Prob_2^*(\R^2)$}
that $m_{\theta,\mu} = \mean(Q_{\theta,\mu}) = \langle m_\mu,\theta\rangle = 0$
and \smash{$\sigma_{\theta,\mu}^2 = \std(Q_{\theta,\mu})^2 = \theta^\top\Sigma_\mu\theta = \theta^\top\theta = 1$}.
For any $\theta \in \Sph^1$.
Thus, by construction \eqref{eq: nrcdt} and Lem.~\ref{lem: orthogonal whitened},
for any $\mu \in \Prob_2^*(\R^2)$,
$X_{\mu}(\cdot,\theta) = Q_{\theta,\Wht[\mu]} \in \Leb^2_\rho(\UniE)$.
Due to this cancellation of the normalization as a post-processing step \eqref{eq: nrcdt},
motivates the naming of the WR-CDTs, cf.~\eqref{eq: wrcdt}.

\section{Complete Characterization of Affine Classes}

Unlike for the NR-CDTs \cite{beckmann2025max,beckmann2026generalizations},
where injectivity does not hold, see Sec.~\ref{sec: non-inj nrcdt},
and the invertibility of the R-CDT holds only on positive densities \cite{kolouri2016radon},
we establish a measure-theoretic characterization
of affine equivalence through orthogonal orbits of whitened quantile profiles
and show the injectivity of the constructed features.

\begin{proposition}
    \label{prop: group action}
    The orthogonal group $\Orth_2(\R)$ acts on \smash{$\Met \coloneqq \Leb^2_{\rho \otimes u_{\Sph^1}}(\UniSphE{1})$} by
    \begin{equation*}
        \smash{\Orth_2(\R) \times \Met \to \Met,
        \quad O \cdot X = X \circ O^\top,
        \quad \text{where} \quad
        X \circ O^\top(t, \theta)= X(t, O^\top \theta).}
    \end{equation*}
    Each $O \in \Orth_2(\R)$ induces a linear isometry of \smash{$\Leb^2_{\rho \otimes u_{\Sph^1}}(\UniSphE{1})$}.
    Therefore, the \emph{orthogonal orbits} are well-defined:
    \begin{equation}
        \label{eq: orth orbit}
        \smash{\Orbit_\mu
        =\{X_\mu \circ O^\top : O\in \Orth_2(\R)\}
        \subset \Leb^2_{u_{\Sph^1}}(\Sph^1; \Leb_\rho^2(\UniE))
        \cong \Met.}
    \end{equation}
\end{proposition}

\begin{proof}
    Denote $\Met \coloneqq \Leb^2_{\rho \otimes u_{\Sph^1}}(\UniSphE{1})$.
    Let $O_1, \Orth_2 \in \Orth_2(\R)$ and $X \in \Met$.
    Then,
    \begin{equation*}
        \smash{(O_1\Orth_2) \cdot X
        = X \circ (O_1\Orth_2)^\top
        = X \circ \Orth_2^\top O_1^\top
        = ( X \circ \Orth_2^\top) \circ O_1^\top
        = O_1 \cdot (\Orth_2 \cdot X).}
    \end{equation*}
    Let $\iota$ be the identity element of $\Orth_2(\R)$.
    Then, we have $\iota \cdot X = X \circ I_2^\top = X$.
    For each $O \in \Orth_2(\R)$, the map
    $(t,\theta)\mapsto(t,O^\top\theta)$
    preserves the product measure
    $\rho\otimes u_{\Sph^1}$,
    since $u_{\Sph^1}$ is invariant under
    orthogonal transformations, cf. \cite[Thm.~2.4.4]{makarov2013real}.
    Consequently, the action is well-defined
    on $\Leb^2$ equivalence classes.
    Furthermore, for every $X \in \Met$,
    \begin{align*}
        \smash{\|O \cdot X\|_{\Lebrhou}^2
        = \int_{\UniSph{1}} \left|X(t,O^\top \theta)\right|^2 \dd(\rho \otimes u_{\Sph^1})(t,\theta)
        = \|X\|_{\Lebrhou}^2,}
    \end{align*}
    where the second equality follows from the invariance of the product measure,
    cf. \cite[Thm.~2.4.4]{makarov2013real}.
    Therefore,
    $O \cdot X \in \Met$, and the action is a linear isometry.
    \qed
\end{proof}

\begin{proposition}
    \label{prop: group_equivalence}
    For $\mu, \nu \in \Prob_2^*(\R^2)$ holds $\Orbit(X_\mu) = \Orbit(X_\nu)$
    if and only if $X_\mu = O \cdot X_\nu$
    for some $O \in \Orth_2(\R)$.
\end{proposition}
\begin{proof}
    Starting from the definition of the orbit $\Orbit(X_\mu)=\{O_1 \cdot X_\mu : O_1 \in \Orth_2(\R)\}$,
    and assume that $X_\nu = O_1 \cdot X_\mu$ for some $O_1 \in \Orth_2(\R)$,
    then, by Prop.~\ref{prop: group action}, holds
    \begin{align*}
        \smash{\Orbit(X_\nu)
        = \{\Orth_2 \cdot (O_1 \cdot X_\mu) \mid \Orth_2 \in \Orth_2(\R)\}
        = \{(\Orth_2O_1) \cdot X_\mu \mid \Orth_2 \in \Orth_2(\R)\}
        = \Orbit(X_\mu),}
    \end{align*}
    because $\Orth_2 \mapsto \Orth_2O_1$ is a bijection of $\Orth_2(\R)$.
    Conversely,
    if $\Orbit(X_\mu)=\Orbit(X_\nu)$,
    then $X_\nu = I \cdot X_\nu \in \Orbit(X_\nu) = \Orbit(X_\mu)$.
    Consequently,
    $X_\nu = O \cdot X_\mu$ for some $O \in \Orth_2(\R)$
    by construction \eqref{eq: orth orbit}.
    \qed
\end{proof}

\begin{lemma}
    \label{lem: orthogonal X}
    If $\mu \in \Prob_2^*(\R^2)$ is whitened and $O \in \Orth_2(\R)$,
    then we have $X_{O_\# \mu}(\cdot,\theta) = X_\mu(\cdot,O^\top\theta)$
    for any $\theta \in \Sph^1$.
    In particular, $X_{O_\# \mu}=O \cdot X_\mu \in \Met = \Leb^2_{\rho \otimes u_{\Sph^1}}(\UniSphE{1})$.
\end{lemma}
\begin{proof}
    Since $O_\#\mu$ is whitened by Lem.~\ref{lem: orthogonal whitened},
    and \eqref{eq: radon-covariance}, we have $\Rad_\theta[O_\#\mu] = \Rad_{O^\top\theta}[\mu]$ for every $\theta \in \Sph^1$.
    Taking CDTs and using \eqref{eq: wrcdt} we get
    \begin{equation*}
        \smash{X_{O_\#\mu}(\cdot,\theta)
        = Q_{\theta, \Wht[O_\#\mu]}
        = Q_{O^\top\theta, \Wht[\mu]}
        = X_\mu(\cdot,O^\top\theta).}
        \tag*{\qed}
    \end{equation*}
\end{proof}

For proving the main result
on the injectivity of the feature extractor,
we provide a variant of the Cramér–Wold theorem \cite{cramer1936distribution}.
Therefore, the Fourier transform of a probability measure $\mu \in \Prob(\R^d)$
is defined as $\Four_\mu (\xi) \coloneqq \int_{\R^d} \exp(-\mathrm i \langle \xi, x\rangle) \dd \mu(x)$.
Hence, by the well-known Fourier-slice theorem \cite{helgason1999radon},
the following relation between the Fourier transform
on the Radon transform $\Rad_\theta[\mu] \in \Prob(\R)$
and the Fourier transform on $\mu$ correspond as follows
\begin{align}
    \label{eq: fourier slice}
    \smash{\Four_{\Rad_\theta[\mu]}(t)
    = \int_{-\infty}^\infty \hspace{-5pt} \exp(- \mathrm i st) \dd \!\!\!\overbracket{\Rad_\theta[\mu]}^{= (\langle \cdot, \theta\rangle)_{\#}\mu}\!\!(s)
    = \int_{\R^d} \hspace{-5pt} \exp(- \mathrm i \langle \theta, x\rangle t) \dd \mu(x)
    = \Four_\mu (t \theta).}
\end{align}

\begin{proposition}
    \label{prop: cramer wold}
    Let $\mu,\nu \in \Prob(\R^d)$ and let
    $\Disk \subset \Sph^{d-1}$ such that $\Sph^{d-1} \setminus \Disk$
    is a $u_{\Sph^{d-1}}$ measure zero set, and
    $\Rad_\theta[\mu] = \Rad_\theta[\nu]$ for all $\theta \in \Disk$.
    Then $\mu = \nu$.
\end{proposition}

\begin{proof}
    By \eqref{eq: fourier slice}
    and similar to \cite[Thm.~1]{boman2009support} %
    we obtain
    \\[2pt]
    \begin{equation*}
        \smash{\Four_\mu(t \theta)
        = \Four_{\Rad_\theta[\mu]}(t)
        = \int_{-\infty}^\infty \exp(-\mathrm i st) \dd \overbracket{\Rad_\theta[\mu]}^{= \Rad_\theta[\nu]}(s)
        = \Four_{\Rad_\theta[\nu]}(t)
        = \Four_\nu(t \theta)}
    \end{equation*}
    for all $\theta \in \Disk$ and $t \in \R$.
    Hence, for any $\xi \in \R^d\setminus\{0\}$,
    we find a sequence $(\theta_n)_{n \in \N} \in \Disk^{\N}$ such that
    $\theta_n \to \xi/\|\xi\|$ as $n \to \infty$
    if and only if
    $\theta_n \|\xi\| \to \xi$ as $n \to \infty$.
    By Lebesgue's dominated convergence theorem,
    $\Four_\mu$ is continuous,
    and consequently, we obtain
    \begin{equation*}
        \smash{\Four_\mu[\xi]
        = \lim_{n \to \infty} \Four_\mu[\theta_n \|\xi\|]
        = \lim_{n \to \infty} \Four_\nu[\theta_n \|\xi\|]
        = \Four_\nu[\xi],
        \quad \forall \xi \in \R^d\setminus\{0\}.}
    \end{equation*}
    For $\xi = 0$, all is equal to one.
    Hence, $\mu = \nu$ by the uniqueness of the Fourier transform \cite[Sec.~26]{billingsley1995probability}
    for finite Borel measures.
    \qed
\end{proof}

\begin{theorem}
    \label{thm: orbit injectivity}
    For $\mu,\nu \in \mathcal{P}_2^*(\R^2)$
    it holds $\Orbit(X_\nu) = \Orbit(X_\mu)$
    if and only if $\nu \in [\mu]$.
\end{theorem}
\begin{proof}
    For ``$\Leftarrow$'':
    let $\nu \in [\mu]$.
    By Prop.~\ref{prop: whitening},
    there exists $O \in \Orth_2(\R)$ such that $\Wht[\nu] = O_{\#}\Wht[\mu]$.
    By Lem.~\ref{lem: orthogonal X},
    $X_\nu = X_{O_{\#}\Wht[\mu]} = O \cdot X_\mu$,
    and by Prop.~\ref{prop: group_equivalence},
    the assertion follows.

    For ``$\Rightarrow$'':
    let $ \Orbit(X_\nu) = \Orbit(X_\mu)$.
    By Prop.~\ref{prop: group_equivalence},
    there exists $O \in \Orth_2(\R)$ such that $X_\nu = O \cdot X_\mu \in \Leb^2_{\rho \otimes u_{\Sph^1}}(\UniSphE{1}) = \Met$.
    We prove at first
    that $\Wht[\nu] = O_\# \Wht[\mu]$,
    to apply then Prop.~\ref{prop: whitened to affine}
    which yields the assertion.
    Therefore,
    by Lem.~\ref{lem: orthogonal whitened}
    both $\Wht[\nu]$ and $O_\# \Wht[\mu]$ are whitened,
    and by Lem.~\ref{lem: orthogonal X}
    it follows that $X_\nu = O \cdot X_\mu \in \Met$.
    To this end,
    by definition and \eqref{eq: wrcdt} it follows that
    \begin{align*}
        \smash{0
        = \|X_{\nu} - O \cdot X_{\mu}\|^2_{\Lebrhou}
        = \int_{\Sph^1} \|Q_{\theta,\Wht[\nu]} - Q_{\theta,O_\# \Wht[\mu]}\|^2_{\Lebrho}
        \dd u_{\Sph^1}(\theta)}
    \end{align*}
    using Tonelli's theorem.
    Using \eqref{eq: sliced wasserstein} yields
    \smash{$0 = \SWass_2^2(\Wht[\nu], O_{\#}\Wht[\mu])$}.
    Thus,
    almost everywhere with respect to $u_{\Sph^1}$
    holds $\Wass_2({\mathcal{R}_\theta[\Wht[\nu]]},{\mathcal{R}_\theta[O_\# \Wht[\mu]]}) = 0$.
    Since $\Wass_2$ is a metric,
    it holds $\Rad_\theta[\Wht[\nu]] = \Rad_\theta[O_\# \Wht[\mu]]$
    for all $\theta \in \Disk$ with $\Disk$ having a $u_{\Sph^1}$ null complement.
    By Prop.~\ref{prop: cramer wold}
    it follows that $\Wht[\nu] = O_\# \Wht[\mu]$.
    This finalizes the proof.
    \qed
\end{proof}

\section{A Metric on Affine Classes}

We equip the space of affine classes $\mathcal A$
from \eqref{eq: aff class fac} with a metric.
Covariance whitening reduces affine transformations to an orthogonal
ambiguity, cf.~Prop.~\ref{prop: whitening} and~\ref{prop: whitened to affine},
which we resolve by optimal sliced Wasserstein alignment.
Retaining the full directional representation avoids the aggregation
ambiguities of $_h$NR-CDTs (§~\ref{sec: non-inj nrcdt})
and yields a distance that vanishes exactly for affinely equivalent
measures.
The induced representation of $\mathcal A$ is thus injective.

\begin{lemma}
\label{lem: density}
    $C_{\mathrm c}(\UniSphE{1},\R)$ is dense in
    $\Leb^2(\UniSphE{1},\mathcal{B}(\UniSphE{1}),\rho\otimes u_{\Sph^1})$
\end{lemma}

\begin{proof}
    The space $\UniSphE{1}$ is a locally compact metric space.
    Since both factors are second-countable,
    $\mathcal B(\UniSphE{1})=\mathcal B(\R)\otimes\mathcal B(\Sph^1)$.
    Thus, $\rho\otimes u_{\Sph^1}$ is a finite Borel measure
    on $\UniSphE{1}$, which is regular by
    \cite[Theorem~2.18]{rudin1987real}.
    The density of continuous functions with compact support in \(\Leb^2\)
    follows from \cite[Theorem~3.14]{rudin1987real}. %
    \qed
\end{proof}

\begin{lemma}
    \label{lem: continuity orth}
    For every $X \in {\Leb^2_{\rho \otimes u_{\Sph^1}}(\UniSphE{1})}$,
    the map $\Orth_2(\R) \to {\Leb^2_{\rho \otimes u_{\Sph^1}}(\UniSphE{1})}: O \mapsto O \cdot X$
    from Prop.~\ref{prop: group action} is continuous.
\end{lemma}
\begin{proof}
    Let $\varepsilon > 0$ and $O_n \to O$ in $\Orth_2(\R)$.
    By Lem.~\ref{lem: density},
    there exists a $Y\in C_{\mathrm c}(\UniSphE{1},\R)$,
    with \smash{$\|X-Y\|_{\raisebox{10pt}{\Lebrhou}} < \varepsilon$}.
    Since $O_n^\top \to O^\top$,
    and for fixed $\theta \in \Sph^1$,
    $\|O_n^\top \theta - O^\top\theta\|_2 \leq \|O_n^\top -O^\top\|_{\op} \|\theta\|_2 = \|O_n^\top -O^\top\|_{\op} \to 0$ as $n \to \infty$,
    pointwise convergence follows
    $(O_n \cdot Y)(t, \theta) \to (O \cdot Y)(t, \theta)$ as $n \to \infty$.
    Since $Y$ is continuous with compact support,
    it is bounded
    $\sup |Y| = c_Y < \infty$.
    Therefore,
    $|(O_n \cdot Y)(t, \theta)-(O \cdot Y)(t, \theta)| \leq 2 c_Y$.
    Since \(\rho\otimes u_{\mathbb S^1}\) is a finite measure,
    it follows from Lebesgue's dominated convergence theorem that
    $\|O_n \cdot Y - O \cdot Y\|_{\Lebrhou} \to 0$ as $n \to \infty$.
    Since the group action is linear,
    it holds
    \begin{equation*}
        \smash{\| O_n \cdot X - O \cdot X \|_{\Lebrhou}
        < \| O_n \cdot Y - O \cdot Y \|_{\Lebrhou}
        + 2\varepsilon.}
    \end{equation*}
    Therefore, the continuity follows by
    \smash{${\lim \sup}_{{n\to\infty}} \| O_n \cdot X - O \cdot X \|_{\raisebox{10pt}{\Lebrhou}} \le 2 \varepsilon$}.
    \hspace{-5pt}
    \qed
\end{proof}

\begin{corollary}
\label{cor: orbit closed}
    Every orbit $\Orbit(X)$ is closed in ${\Leb^2_{\rho \otimes u_{\Sph^1}}(\UniSphE{1})}$.
\end{corollary}
\begin{proof}
    Since $\Orth_2(\R)$ is compact,
    and the map \smash{$\Orth_2(\R) \to {\Leb^2_{\rho \otimes u_{\Sph^1}}(\UniSphE{1})}: O \mapsto O \cdot X$}
    is continuous by Lem.~\ref{lem: continuity orth},
    $\Orbit(X)$ is compact,
    and thus closed in ${\Leb^2_{\rho \otimes u_{\Sph^1}}(\UniSphE{1})}$.
    \qed
\end{proof}

Let $(\Met,d)$ be a nonempty metric space.
For nonempty $A, B \subset \Met$,
we define the distance from $A$ to $B$ in $\Met$ \cite[§~3.3, pp.~76--77]{burago2001course}
by $\dist(A,B) \coloneqq \inf\{d(X,Y) \mid X \in A,\ Y \in B\}$.
Let a group $\Gr$ act on $\Met$ by isometry,
denote the orbit of $X \in \Met$ by $\Gr(X)$,
and let $\Met/\Gr \coloneqq \{\Gr(X) : X \in \Met\}$ be the set of orbits.
The \emph{orbit space distance function} is defined \cite[§~3.3, pp.~76--77]{burago2001course} by
\begin{equation}
    \label{eq: orbit distance}
    \smash{d_{\Gr} : \Met / \Gr \times \Met / \Gr \to \R, \qquad d_{\Gr}(\Gr(X), \Gr(Y))  \coloneqq \dist (\Gr(X), \Gr(Y)).}
\end{equation}

\begin{theorem}
    \label{thm: metric orth}
    For \smash{$\Met = \Leb^2_{\rho\otimes u_{\Sph^1}}(\UniSphE{1})$}
    let $d_{\Orth_2(\R)}$ be the orbit space distance function \eqref{eq: orbit distance}
    with \smash{$d(\cdot_1,\cdot_2) = \|\cdot_1 - \cdot_2\|_{\raisebox{10pt}{\Lebrhou}}$},
    $\Gr = \Orth_2(\R)$, and $\Gr(X)=\Orbit(X)$.
    Then $d_{\Orth_2(\R)}$ is a metric on $\Leb^2_{\rho\otimes u_{\Sph^1}}(\UniSphE{1})/\Orth_2(\R)$,
    and for all $X,Y \in \Leb^2_{\rho\otimes u_{\Sph^1}}(\UniSphE{1})$,
    \begin{equation*}
        \smash{d_{\Orth_2(\R)}\bigl(\Orbit(X),\Orbit(Y)\bigr)
        = \dist\bigl(\{X\}, \Orbit(Y)\bigr)
        = \inf_{O\in \Orth_2(\R)}\,\|X - O \cdot Y\|_{\Lebrhou}.}
    \end{equation*}
\end{theorem}
\begin{proof}
    By Prop.~\ref{prop: group action},
    $\Orth_2(\R)$ acts on \smash{$\bigl(\Leb^2_{\rho\otimes u_{\Sph^1}}(\UniSphE{1}),\ \|\cdot\|_{\raisebox{10pt}{\Lebrhou}}\bigr)$}
    by isometries, so $d(O_1 \cdot X, O_2 \cdot Y) = d(X,(O_1^{-1}O_2) \cdot Y)$ for all $O_1,O_2 \in \Orth_2(\R)$.
    Using $\{O_1^{-1}O_2 : O_1,O_2 \in \Orth_2(\R)\} = \Orth_2(\R)$
    and the infimum over $O_1,O_2 \in \Orth_2(\R)$
    gives $\dist(\Orbit(X),\Orbit(Y)) = \dist(\{X\},\Orbit(Y))$.
    By Lem.~\ref{lem: continuity orth}
    and compactness of \(O_2(\mathbb R)\), the infimum is attained
    By Cor.~\ref{cor: orbit closed} every orbit is closed,
    hence $d_{\Orth_2(\R)}$ is a metric, by \cite[Thm.~6.6.1, p.~235]{ratcliffe2006manifolds}.
    \qed
\end{proof}

\begin{corollary}
\label{thm: metric aff}
    The function $d_\mathcal{A}: \mathcal{A} \times \mathcal{A} \to \R$,
    $([\mu],[\nu]) \mapsto d_ {\Orth_2(\R)}\bigl(\Orbit(X_\mu),\Orbit(X_\nu)\bigr)$
    is well defined
    and a metric on the set $\mathcal{A}$ of affine classes \eqref{eq: aff class fac}.
\end{corollary}

\begin{corollary}
    \label{margins}
    Let $\mu_1,\ldots,\mu_c \in \Prob_2^*(\R^2)$ for $c \ge 2$
    lie in pairwise distinct affine classes, i.e. $\mu_i \not\sim \mu_j, i \neq j$.
    Then, it holds $\gamma \coloneqq \min_{i \neq j} d_\mathcal{A}([\mu_i],[\mu_j]) > 0$
    and
    \\[-15pt]
    \begin{enumerate}
        \item[i)] If $\nu\in[\mu_i]$,
        then $d_ {\Orth_2(\R)}(\Orbit(X_\nu),\Orbit(X_{\mu_i})) = 0$
        and $d_ {\Orth_2(\R)}(\Orbit(X_\nu),\Orbit(X_{\mu_j})) \geq \gamma$,
        \\[-10pt]
        \item[ii)] Let $\widetilde{X}\in \Met$
        with
        $d_{\Orth_2(\R)}(\Orbit(\widetilde{X}),\Orbit(X_{\mu_i})) < \tfrac{\gamma}{2}$,
        then
        $d_{\Orth_2(\R)}(\Orbit(\widetilde{X}),\Orbit(X_{\mu_j})) > \tfrac{\gamma}{2}$.
    \end{enumerate}
\end{corollary}

\section{Numerical Experiments}

All descriptors use a Python implementation based on
\cite{beckmann2026generalizations}\footnote{Original code is publicly available at \href{https://github.com/DrBeckmann/NR-CDT}{GitHub}},
with a discrete Radon transform followed by the CDT,
using $N$ angles,
the same radii for the datasets as in \cite{beckmann2025max,beckmann2026robustness},
$Q=64$ quantile levels
and $\rho = u_{[0,1]}$,
so integration in the reference variable reduces to integration over \([0,1]\).
For WR-CDT, we use the directions $\theta_\mu(\theta_k)$
of Cor.~\ref{cor: no-resampling} for equispaced $\theta_k\in\Sph^1$,
correct the CDT's half-step quantile offset,
and normalize by exact projection means and standard deviations
(Lem.~\ref{lem: mean,std}) to obtain $X[\mu]$.
All methods share the data, splits, queries, $N$, and $Q$,
differing only in directions and profile comparison\footnote{Code and experiments are publicly available at \href{https://github.com/christos-sonis/WR-CDT}{GitHub}}.

\paragraph{Discrete orbit distance.}
We approximate \eqref{eq: orbit distance} by minimizing the RMS
discrepancy over circular shifts and angular reversals of the sampled
profiles, representing grid rotations and reflections.
For $X=X[\mu]$ and $Y=X[\nu]$, the squared discrepancy at shift $s$ is
$(\|X\|_{\mathrm F}^2+\|Y\|_{\mathrm F}^2-2C_s)/(QN)$, where
$C_s=\sum_{q=1}^Q\sum_{j=0}^{N-1}X_{q,j}Y_{q,(j-s)\bmod N}$.
We obtain all $C_s$ by inverse DFT of the summed rowwise
cross-spectra \cite{cooley1969finite}, repeating the calculation
for the reflected profile.
The FFT preserves the exact grid minimum while reducing
pairwise alignment from $\mathcal{O}(QN^2)$
to $\mathcal{O}(QN\log N)$ \cite{cooley1965algorithm}.

\begin{figure}[t!]
    \includegraphics[width=\linewidth, trim={10pt 42pt 0pt 12pt}, clip]{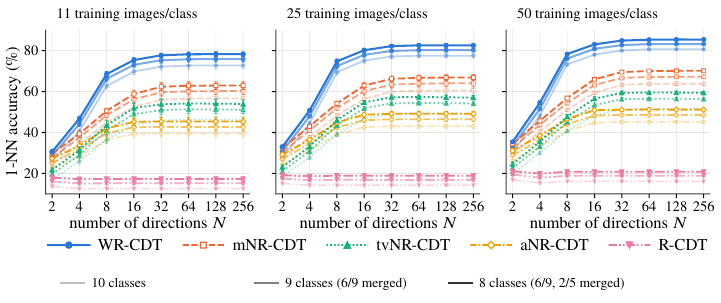}
    \includegraphics[width=\linewidth, trim={10pt 0pt 5pt 0pt}, clip]{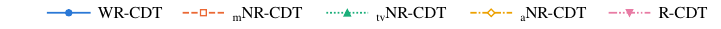} \\[-5pt]
    \includegraphics[width=\linewidth, trim={-30pt 5pt -10pt 133pt}, clip]{numeric/fig_linmnist_columnwise_overlay.pdf}
    \\[-23pt]
    \caption{LinMNIST $1$-NN accuracy (mean$\pm$std, \%)
    versus number of angles $N$ with $11$, $25$, and $50$
    training images per class (left to right)
    with ten classes and merged classes $(6/9), (2/5)$.
    }
    \label{fig: linmnist}
    \vspace{-10pt}
\end{figure}

\paragraph{LinMNIST classification.}
We evaluate $1$-nearest-neighbor ($1$-NN) classification
on all ten LinMNIST \cite{BeckmannHeilenkoetter2024} digit classes,
with $500$ randomly affinely transformed images per class.
For each of $20$ random splits at each training size,
we select $11$, $25$, or $50$ training images per class
and classify the remaining images.
Each test image receives the label of the training image at minimum distance.
We use the discrete orbit distance for WR-CDT and Euclidean distance for R-CDT
and the $_{\mathrm{m}}$NR-CDT \cite{beckmann2025max}, $_{\mathrm{tv}}$NR-CDT \cite{beckmann2026generalizations}, and $_{\mathrm{a}}$NR-CDT \cite{beckmann2026robustness} descriptors.
Fig.~\ref{fig: linmnist} shows mean accuracy and standard deviation across splits.
For the best setting of $50$ training images per class,
WR-CDT reaches $73.28\%$ using eight angles
and saturates at $80\%$ with $32$ or more angles,
whereas the highest baseline accuracy across the tested grids is $63.82\%$.
Merging, firstly, the classes $6$ and $9$, being ideally rotated versions,
and secondly, the classes $2$ and $5$, being reflected versions,
increases the accuracy to $83.21\%$ ($+3.12\%$) and $85.40\%$ ($+5.84\%$), respectively.

\paragraph{Chinese character classification.}
\begin{wrapfigure}{r}{0.35\linewidth}
    \vspace{-0pt}
    \includegraphics[width=\linewidth, trim={0pt 42pt 0pt 0pt}, clip]{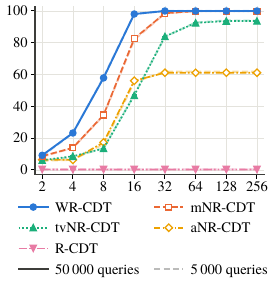}
    \includegraphics[width=\linewidth, trim={0pt 0pt 0pt 125pt}, clip]{numeric/fig_chinese_overlay_nolabels.pdf}
    \\[-25pt]
    \caption{Chinese character accuracy (mean$\pm$std, \%)
    versus number of angles for $5\,000$ and $50\,000$ queries.
    Labels as in Fig.~\ref{fig: linmnist}.
    }
    \vspace{-10pt}
    \label{fig: chinese}
\end{wrapfigure}
We evaluate nearest-template classification
with one reference image for each of $1\,000$ Chinese characters ~\cite{Bliem2020Hanzi}
and $5$ or $50$ random affine copies of each as queries
(in total $5\,000$ or $50\,000$). %
We assign each query to its closest template
using the same representations and distances as for LinMNIST.
Fig.~\ref{fig: chinese} reports accuracy with mean and standard deviation
using $N = 256$ angles for both query-sets
and shows the effect of varying the number of angles for the larger set.
The accuracies for the two query-set sizes are similar.
Additionally, the WR-CDT needs approximately half as many projection angles,
to achieve the same accuracies as by the $_\mathrm{m}$NR-CDT.

\begin{figure}[htb] %
    \begin{minipage}[c]{0.4\linewidth}
        \vspace{-10pt}
        \captionof{table}{Recovery errors by $\kappa$.}
        \vspace{1pt}
        \label{tab3}
        \resizebox{\linewidth}{!}{%
            \begin{tabular}{l@{\quad}r@{\quad}r@{\quad}r@{\quad}r}
            \toprule
            $\kappa$ %
            & \scriptsize{rel.\ RMS error (\%)}
            & \scriptsize{rel.\ matrix error (\%)}
            & \scriptsize{rel.\ RMS} \\
            & \scriptsize{med. / 95th pct.}
            & \scriptsize{med. / 95th pct.} & \scriptsize{$< 2\%$} \\
            \midrule
            1 %
            & 0.23 / 0.55 & 0.18 / 0.48 & 100\% \\
            2 %
            & 0.27 / 0.66 & 0.20 / 0.54 & 99.9\% \\
            4 %
            & 0.31 / 1.03 & 0.22 / 0.76 & 99.2\% \\
            8 %
            & 0.33 / 1.61 & 0.24 / 1.14 & 97.3\% \\
            \bottomrule
            \end{tabular}
        }
    \end{minipage}\hfill%
    \begin{minipage}[c]{0.58\linewidth}
        \centering
        \includegraphics[width=\linewidth, trim={0pt 2pt 0pt 9pt}, clip]{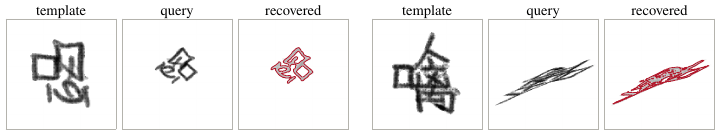} \\[-12pt]
        \captionof{figure}{Affine recovery with reflection at
        $\kappa=1$ (left, relative RMS error $0.28\%$)
        and at $\kappa=8$ (right, median error $0.33\%$).
        Each block shows the template, bilinearly resampled query,
        and query with the estimated contour.}
        \label{fig: recovery}
    \end{minipage}
\end{figure}

\vspace{-1pt}
\paragraph{Affine map recovery.}
We recover $x\mapsto Ax+y$ from $5\,000$ rasterized affine copies
of Chinese character templates without near-symmetries,
with condition number $\kappa$.
Using $N = 256$ angles, we refine the grid minimizer
of \eqref{eq: orbit distance} between neighboring angles
to obtain the orthogonal alignment $\widehat O$ from template to query.
Undoing whitening gives
\smash{$\widehat A=\Sigma_\nu^{1/2}\widehat O\Sigma_\mu^{-1/2}$}
and \smash{$\widehat y=m_\nu-\widehat A m_\mu$},
cf.~Prop.~\ref{prop: whitened to affine}.
Tab.~\ref{tab3} reports the RMS displacement of
$\widehat A x+\widehat y$ from $Ax+y$ over $x\sim\mu$,
relative to the query's RMS radius, and the relative matrix error
$\|\widehat A-A\|_{\mathrm F}/\|A\|_{\mathrm F}$.
Fig.~\ref{fig: recovery} illustrates two recovered alignments.
At $\kappa=8$, the median relative RMS error is only $0.33\%$,
with $97.3\%$ of recovered maps below $2\%$ error.
\\[-25pt]

\vspace{-25pt}
\section{Conclusion}

We introduced an affine-invariant metric based on an injective
representation of affine classes.
Whitening reduces affine alignment to orthogonal alignment
while preserving the $\Leb^2$ geometry of the sliced Wasserstein distance.
The minimum
$d_\mathcal{A}([\mu],[\nu])
= \min_{O\in\Orth_2(\R)}
\SWass_2\bigl(\Wht[\mu],O_\#\Wht[\nu]\bigr)$
is attained, making nearest-orbit classification well posed.
The underlying theoretical construction naturally generalizes to higher dimensions
Experiments demonstrate improved nearest-neighbor classification,
comparable accuracy with about half as many angles
in Chinese character classification, and accurate affine map recovery.
These results motivate investigating whether the metric
can reduce reliance on affine data augmentation in machine learning approaches.
\\[-15pt]

\bibliographystyle{splncs04}
\bibliography{bib-doi}

@book{villani2003topics,
  author = {Villani, C{\'e}dric},
  title = {Topics in Optimal Transportation},
  publisher = {Amer. Math. Soc.},
  year = {2003},
  doi = {10.1090/gsm/058}
}

@article{peyre2019computational,
  author = {Peyr{\'e}, Gabriel and Cuturi, Marco},
  title = {Computational Optimal Transport},
  journal = {Found. Trends Mach. Learn.},
  year = {2019},
  volume = {11},
  pages = {355--607},
  doi = {10.1561/2200000073}
}

@article{park2018cdt,
  author = {Park, Se Rim and Kolouri, Soheil and Kundu, Shinjini and Rohde, Gustavo K.},
  title = {The Cumulative Distribution Transform and Linear Pattern Classification},
  journal = {Appl. Comput. Harmon. Anal.},
  year = {2018},
  volume = {45},
  pages = {616--641},
  doi = {10.1016/j.acha.2017.02.002}
}

@article{kolouri2016radon,
  author = {Kolouri, Soheil and Park, Se Rim and Rohde, Gustavo K.},
  title = {The {Radon} Cumulative Distribution Transform and Its Application to Image Classification},
  journal = {IEEE Trans. Image Process.},
  year = {2016},
  volume = {25},
  pages = {920--934},
  doi = {10.1109/TIP.2015.2509419}
}

@article{bonneel2015sliced,
  author = {Bonneel, Nicolas and Rabin, Julien and Peyr{\'e}, Gabriel and Pfister, Hanspeter},
  title = {Sliced and {Radon} {Wasserstein} Barycenters of Measures},
  journal = {J. Math. Imaging Vis.},
  year = {2015},
  volume = {51},
  pages = {22--45},
  doi = {10.1007/s10851-014-0506-3}
}

@inproceedings{beckmann2025max,
  author = {Beckmann, Matthias and Beinert, Robert and Bresch, Jonas},
  title = {Max-Normalized {Radon} Cumulative Distribution Transform for Limited Data Classification},
  booktitle = {{SSVM} 2025},
  series = {LNCS},
  volume = {15667},
  year = {2025},
  pages = {241--254},
  publisher = {Springer},
  doi = {10.1007/978-3-031-92366-1_19}
}

@article{beckmann2026robustness,
  author = {Beckmann, Matthias and Beinert, Robert and Bresch, Jonas},
  title = {Normalized {Radon} Cumulative Distribution Transforms for Invariance and Robustness in Optimal Transport Based Image Classification},
  journal = {SIAM J. Imaging Sci.},
  year = {2026},
  volume = {19},
  pages = {677--709},
  doi = {10.1137/25M1767741}
}

@article{beckmann2026generalizations,
  author = {Beckmann, Matthias and Beinert, Robert and Bresch, Jonas},
  title = {Generalizations of the Normalized {Radon} Cumulative Distribution Transform for Limited Data Recognition},
  journal = {J. Math. Imaging Vis.},
  year = {2026},
  volume = {68},
  pages = {45},
  doi = {10.1007/s10851-026-01327-4}
}

@article{lan2013whitening,
  author = {Lan, Rushi and Yang, Jianwei and Jiang, Yong and Fyfe, Colin and Song, Zhan},
  title = {Whitening Central Projection Descriptor for Affine-Invariant Shape Description},
  journal = {IET Image Process.},
  year = {2013},
  volume = {7},
  pages = {81--91},
  doi = {10.1049/iet-ipr.2012.0094}
}

@book{santambrogio2015optimal,
  author = {Santambrogio, Filippo},
  title = {Optimal Transport for Applied Mathematicians: Calculus of Variations, {PDEs}, and Modeling},
  publisher = {Birkh{\"a}user},
  year = {2015},
  doi = {10.1007/978-3-319-20828-2}
}

@article{cramer1936distribution,
  author = {Cramér, Harald and Wold, Herman},
  title = {Some Theorems on Distribution Functions},
  journal = {J. Lond. Math. Soc.},
  volume = {s1-11},
  pages = {290--294},
  year = {1936},
  doi = {10.1112/jlms/s1-11.4.290}
}

@article{boman2009support,
  author = {Boman, Jan and Lindskog, Filip},
  title = {Support Theorems for the {Radon} Transform and {Cram{\'e}r--Wold} Theorems},
  journal = {J. Theor. Probab.},
  year = {2009},
  volume = {22},
  pages = {683--710},
  doi = {10.1007/s10959-008-0151-0}
}

@book{helgason1999radon,
  author = {Helgason, Sigurdur},
  title = {The {Radon} Transform},
  edition = {2},
  publisher = {Birkh{\"a}user},
  year = {1999},
  doi = {10.1007/978-1-4757-1463-0}
}

@article{cooley1969finite,
  author = {Cooley, James W. and Lewis, Peter A. W. and Welch, Peter D.},
  title = {The Finite {Fourier} Transform},
  journal = {IEEE Trans. Audio Electroacoust.},
  year = {1969},
  volume = {17},
  pages = {77--85},
  doi = {10.1109/TAU.1969.1162036}
}

@article{cooley1965algorithm,
  author = {Cooley, James W. and Tukey, John W.},
  title = {An Algorithm for the Machine Calculation of Complex {Fourier} Series},
  journal = {Math. Comput.},
  year = {1965},
  volume = {19},
  pages = {297--301},
  doi = {10.1090/S0025-5718-1965-0178586-1}
}

@book{billingsley1995probability,
  author = {Billingsley, Patrick},
  title = {Probability and Measure},
  edition = {3},
  year = {1995},
  publisher = {Wiley}
}

@book{ratcliffe2006manifolds,
  author = {John G. Ratcliffe},
  title = {Foundations of Hyperbolic Manifolds},
  edition = {2},
  publisher = {Springer},
  year = {2006},
  doi = {10.1007/978-0-387-47322-2}
}

@book{rudin1987real,
  author = {Walter Rudin},
  title = {Real and Complex Analysis},
  edition = {3},
  publisher = {McGraw-Hill},
  year = {1987}
}

@book{burago2001course,
  author = {Burago, Dmitri and Burago, Yuri and Ivanov, Sergei},
  title = {A Course in Metric Geometry},
  publisher = {Amer. Math. Soc.},
  year = {2001},
  doi = {10.1090/gsm/033}
}

@book{makarov2013real,
  author = {Makarov, Boris and Podkorytov, Anatolii},
  title = {Real Analysis: Measures, Integrals and Applications},
  publisher = {Springer},
  year = {2013},
  doi = {10.1007/978-1-4471-5122-7}
}

@article{BeckmannHeilenkoetter2024,
  author = {Beckmann, Matthias and Heilenk{\"o}tter, Nick},
  title = {Equivariant Neural Networks for Indirect Measurements},
  journal = {SIAM J. Math. Data Sci.},
  volume = {6},
  pages = {579--601},
  year = {2024},
  doi = {10.1137/23M1582862}
}

@misc{Bliem2020Hanzi,
  author = {Bliem, Pascal},
  title = {Handwritten Chinese Character ({Hanzi}) Datasets},
  year = {2020},
  howpublished = {Kaggle},
  url = {https://www.kaggle.com/datasets/pascalbliem/handwritten-chinese-character-hanzi-datasets},
  note = {Accessed 2025-03-26}
}

\end{document}